\documentclass[review]{elsarticle}

\usepackage{enumitem}
\usepackage{amsmath}
\usepackage{amssymb, tikz, float}
\usepackage{verbatim}

\usepackage{lineno,hyperref}
\modulolinenumbers[5]

\journal{Discrete Mathematics}

\newtheorem{lem}{Lemma}
\newtheorem{thm}{Theorem}

\newdefinition{df}{Definition}
\newdefinition{rmk}{Remark}

\newproof{pf}{Proof}
\newproof{pf2}{Proof of Theorem \ref{major}}

\begin{document}

\begin{frontmatter}

\title{Ricci-flat graphs with girth four}
\tnotetext[t1]{This work was supported by the National Natural Science Foundation of China (No. 11201496 and No. 11601093).}



\author[gut]{Weihua He}

\author[sysu]{Jun Luo}

\author[sysu]{Chao Yang\corref{cor1}}

\author[sysu]{Wei Yuan}

\cortext[cor1]{Corresponding author: yangchao0710@gmail.com, yangch8@mail.sysu.edu.cn.}

\address[gut]{Department of Applied Mathematics,
Guangdong University of Technology, Guangzhou, China}

\address[sysu]{School of Mathematics, Sun Yat-Sen University, Guangzhou, China}

\begin{abstract}
Lin-Lu-Yau introduced an interesting notion of Ricci curvature for graphs and obtained a complete characterization for all Ricci-flat graphs with girth at least five \cite{Lin2014671}.
In this paper, we propose a concrete approach to construct an infinite family of distinct  Ricci-flat graphs of girth four with edge-disjoint $4$-cycles and completely characterize all Ricci-flat graphs of girth four with vertex-disjoint $4$-cycles.
\end{abstract}


\end{frontmatter}

\linenumbers

\section{Introduction}
A manifold is Ricci-flat if the Ricci curvature vanishes everywhere. Calabi-Yau manifolds are a special type of Ricci-flat manifolds, which provide a potential model to describe the physical world \cite{Yau2010}. The study of Ricci curvature on manifolds has inspired several attempts to bring the concept of Ricci curvature to graphs. Ollivier introduced the Ricci curvature of Markov chains on metric spaces, including graphs \cite{Ollivier2009810}. By modifying Ollivier's definition, Lin-Lu-Yau proposed a slightly different definition for Ricci curvature of graphs \cite{Lin2011}. These new concepts have received considerable discussions. See for example \cite{Bhattacharya201523,Bauer20172033}. These new notions even find applications within combinatorics and computer science. Among others, we refer to \cite{Cho2013916} for the relation between Ollivier's Ricci curvature and the coloring of graphs and  \cite{Ni2015} for the employment of Ricci curvature in understanding the Internet topology.

This paper considers the Ricci curvature in the sense of Lin-Lu-Yau. We are especially interested in \textit{Ricci-flat} graphs, whose Ricci curvature vanishes on every edge. A very recent pioneering work by Lin-Lu-Yau \cite{Lin2014671} completely characterizes all Ricci-flat graphs with girth at least five.

\begin{thm}
A Ricci-flat graph with girth at least five is isomorphic to: (1) the infinite path, (2) a cycle of length at least six, (3) the dodecahedral graph, (4) the half-dodecahedral graph, or (5) the Petersen graph.
\end{thm}

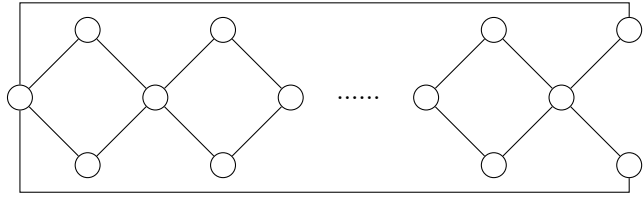
\begin{figure}[h]
\begin{center}
\begin{tikzpicture}[scale=0.18]

\node (a1) at (5,10) [circle,draw] {};
\node (a2) at (0,5) [circle,draw] {};
\node (a3) at (10,5) [circle,draw] {};
\node (a4) at (5,0) [circle,draw] {};

\node (b1) at (15,10) [circle,draw] {};

\node (b3) at (20,5) [circle,draw] {};
\node (b4) at (15,0) [circle,draw] {};

\node at (25,5) {......};

\node (x1) at (35,10) [circle,draw] {};
\node (x2) at (30,5) [circle,draw] {};
\node (x3) at (40,5) [circle,draw] {};
\node (x4) at (35,0) [circle,draw] {};

\node (y1) at (45,10) [circle,draw] {};
\node (y4) at (45,0) [circle,draw] {};

\draw (a1) -- (a2);
\draw (a1) -- (a3);
\draw (a4) -- (a2);
\draw (a4) -- (a3);

\draw (b1) -- (a3);
\draw (b1) -- (b3);
\draw (b4) -- (a3);
\draw (b4) -- (b3);

\draw (x1) -- (x2);
\draw (x1) -- (x3);
\draw (x4) -- (x2);
\draw (x4) -- (x3);
\draw (y1) -- (x3);
\draw (y4) -- (x3);
\draw (y1) -- (45,12) --(0,12)-- (a2);
\draw (y4) --(45,-2) --(0,-2)-- (a2);

\end{tikzpicture}
\end{center}
\caption{A family of Ricci-flat graphs with edge-disjoint $4$-cycles} \label{fig_edge}
\end{figure}

The authors of \cite{Lin2014671} also gave infinitely many examples of Ricci-flat graphs with girth four. We note that, in all their examples, the $4$-cycles have common edges. Even if $4$-cycles having common edges are not allowed, we can still construct infinitely many Ricci-flat graphs with edge-disjoint $4$-cycles, see Figure \ref{fig_edge}.

Then, the remaining Ricci-flat graphs with girth four are those in which every two $4$-cycles are vertex-disjoint, {\em i.e.} having no common vertices. For those graphs, we obtain the following simple characterization, which is the main result of this paper.

\begin{thm}\label{major}
A Ricci-flat graph with girth four such that no vertex is shared by two $4$-cycles  is isomorphic to one of the following two graphs.
\begin{figure}[h]
\begin{center}
\begin{tikzpicture}[scale=0.15]

\node (v0) at (30,30) [circle,draw] {};
\node (v1) at (30,0) [circle,draw] {};
\node (v2) at (0,0) [circle,draw] {};
\node (v3) at (0,30) [circle,draw] {};
\node (v4) at (10,15) [circle,draw] {};
\node (v5) at (20,15) [circle,draw] {};
\node (v6) at (5,22) [circle,draw] {};
\node (v7) at (10,22) [circle,draw] {};
\node (v8) at (5,8) [circle,draw] {};
\node (v9) at (10,8) [circle,draw] {};
\node (v10) at (20,8) [circle,draw] {};
\node (v11) at (25,8) [circle,draw] {};
\node (v12) at (20,22) [circle,draw] {};
\node (v13) at (25,22) [circle,draw] {};

\draw [ultra thick,blue] (v1) -- (v0);
\draw [ultra thick,blue]  (v2) -- (v1);
\draw [ultra thick,blue]  (v3) -- (v0);
\draw [ultra thick,blue]  (v3) -- (v2);
\draw (v6) -- (v3);
\draw (v6) -- (v4);
\draw (v7) -- (v3);
\draw (v7) -- (v5);
\draw (v8) -- (v2);
\draw (v8) -- (v4);
\draw (v9) -- (v2);
\draw (v9) -- (v5);
\draw (v10) -- (v1);
\draw (v10) -- (v4);
\draw (v11) -- (v1);
\draw (v11) -- (v5);
\draw (v12) -- (v0);
\draw (v12) -- (v4);
\draw (v13) -- (v0);
\draw (v13) -- (v5);


\node (n0) at (0+40,30) [circle,draw] {};
\node (n1) at (30+40,30) [circle,draw] {};
\node (n2) at (30+40,0) [circle,draw] {};
\node (n3) at (0+40,0) [circle,draw] {};
\node (n4) at (8+40,10) [circle,draw] {};
\node (n5) at (15+40,20) [circle,draw] {};
\node (n6) at (15+40,10) [circle,draw] {};
\node (n7) at (22+40,20) [circle,draw] {};
\node (n8) at (8+40,20) [circle,draw] {};
\node (n9) at (8+40,25) [circle,draw] {};
\node (n10) at (22+40,5) [circle,draw] {};
\node (n11) at (22+40,10) [circle,draw] {};

\draw [ultra thick,blue]  (n1) -- (n0);
\draw [ultra thick,blue]  (n2) -- (n1);
\draw [ultra thick,blue]  (n3) -- (n0);
\draw [ultra thick,blue]  (n3) -- (n2);
\draw (n4) -- (n3);
\draw [ultra thick,blue]   (n5) -- (n4);
\draw [ultra thick,blue]   (n6) -- (n4);
\draw (n7) -- (n1);
\draw [ultra thick,blue]  (n7) -- (n5);
\draw [ultra thick,blue]  (n7) -- (n6);
\draw (n8) -- (n0);
\draw (n8) -- (n6);
\draw (n9) -- (n0);
\draw (n9) -- (n5);
\draw (n10) -- (n2);
\draw (n10) -- (n6);
\draw (n11) -- (n2);
\draw (n11) -- (n5);

\end{tikzpicture}
\end{center}
\caption{The graphs $R_1$ (left) and $R_2$ (right)} \label{fig_main1b}
\end{figure}
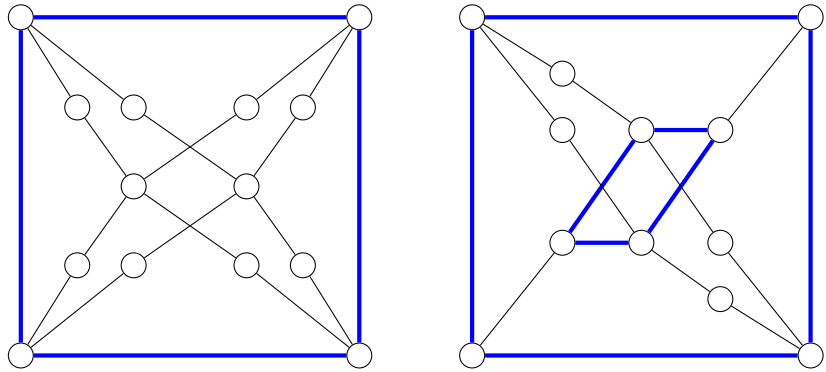
\end{thm}

The rest of this paper is arranged as follows. Section \ref{sec_def} recalls the notion of Ricci curvature on graphs. Section \ref{sec_local} reviews the local structures of Ricci-flat graphs. And Section \ref{sec_main} presents the proof of Theorem \ref{major}.

\section{Preliminaries}\label{sec_def}

We follow Lin-Lu-Yau for the definition of Ricci curvature \cite{Lin2011,Lin2014671}. Let $G$ be a simple undirected graph with vertex set $V$ and edge set $E$. For $x,y\in V$, let $N(x)$ be the set of neighbors of $x$, $d_x=|N(x)|$ be the degree of vertex $x$, and $d(x,y)$ be the distance between $x$ and $y$ in $G$.

A probability distribution is a function $m: V \rightarrow [0,1]$ with $\sum_{x \in V} m(x) = 1$. To define Ricci curvature for each edge of the graph, we only consider distributions $m_x^\alpha$ in the following form,
\begin{equation*}
    m_x^\alpha(v) = \begin{cases}
               \alpha,               & v = x;\\
               \frac{1-\alpha}{d_x},               & x \in N(x);\\
               0, & \text{otherwise},
           \end{cases}
\end{equation*}
where $\alpha \in [0,1]$ and  $x \in V$.

Let $xy\in E$, and let $m_x^\alpha$ and $m_y^\alpha$ be two distributions. A transportation problem between the two distributions can be stated as a linear programming problem. That is, to find the minimum transportation distance
$$\min \sum_{u,v\in V} d(u,v)X_{uv},$$
subject to the constraints
\begin{equation*}
 \begin{cases}
               \sum_{v\in V}  X_{uv} = m_x^\alpha (u),             & u\in V;\\
               \sum_{u\in V}  X_{uv} = m_y^\alpha (v),             & v\in V;\\
               X_{uv} \geqslant 0,
 \end{cases}
\end{equation*}
where the variable $X_{uv}$ denotes the amount transfered from vertex $u$ to vertex $v$.

Define the transportation distance between $m_x^\alpha$ and $m_y^\alpha$ to be optimal solution to the above linear programming problem, namely
$$W(m_x^\alpha, m_y^\alpha)=\min \sum_{u,v\in V} d(u,v)X_{uv}.$$

Ricci curvature is definable on any unordered pair of vertices $x$ and $y$, but for our purpose, we only need the case that $x$ and $y$ are adjacent. For any edge $xy \in E$, the Ricci curvature $\kappa(x,y)$ is defined to be
$$\kappa(x,y)=\lim_{\alpha \rightarrow 1}\frac{1-W(m_x^\alpha, m_y^\alpha)}{1-\alpha}.$$
Recall that a graph $G$ is Ricci-flat if $\kappa(x,y)=0$ for all edges $xy\in E$.

A function $f$ over the vertex set $V$ of $G$ is said to be $c$-Lipschitz if $|f(u)-f(v)|\leqslant c\cdot d(u,v)$ for all $u,v \in V$. By the theory of linear programming, the dual problem of the above defined transportation problem between $m_x^\alpha$ and $m_y^\alpha$ is to find the maximum value
$$\max \sum_{u \in V} f(u)(m_x^\alpha(u)-m_y^\alpha(u))$$
subject to
$$|f(u)-f(v)|\leqslant d(u,v), u,v \in V.$$
In other words, the maximum is taken over all $1$-Lipschitz function $f$. Because the optimal solution of a linear programming problem is equal to that of its dual problem, we have
$$W(m_x^\alpha, m_y^\alpha)=\max \sum_{u \in V} f(u)(m_x^\alpha(u)-m_y^\alpha(u)).$$
Thus, we have the following lemma.

\begin{lem}[\cite{Lin2014671}]\label{lem_dual} Let $f$ be any $1$-Lipschitz function, then
$$W(m_x^\alpha, m_y^\alpha)\geqslant \sum_{u \in V} f(u)(m_x^\alpha(u)-m_y^\alpha(u)).$$
\end{lem}



\section{Local Structures}\label{sec_local}

Before our discussion on the local structure of  Ricci-flat graphs of girth $4$, we recall a lemma from \cite{Lin2014671}.

\begin{lem}[\cite{Lin2014671}]\label{lem_local5}
Suppose that an edge $xy$ in a graph $G$ is not in any $3$-cycles or $4$-cycles, and assume $d_x \leqslant d_y$, then one of the following statements holds.
\begin{enumerate}
\item $d_x=d_y=2$, and $xy$ is not in any $5$-cycle.
\item $d_x=d_y=3$, and $xy$ is shared by two $5$-cycles.
\item $d_x=2, d_y=3$. Let $x_1$ be the other neighbor of $x$ besides $y$, and let $y_1$ and $y_2$ be the two neighbors of $y$ besides $x$, then $\{d(x_1,y_1),d(x_1,y_2)\}=\{2,3\}$.
\item $d_x=2,d_y=4$. Let $x_1$ be the other neighbor of $x$ besides $y$, and let $y_1$, $y_2$ and $y_3$ be the three neighbors of $y$ besides $x$, then at least two of $y_1, y_2, y_3$ have distance $2$ from $x_1$.
\end{enumerate}
\end{lem}

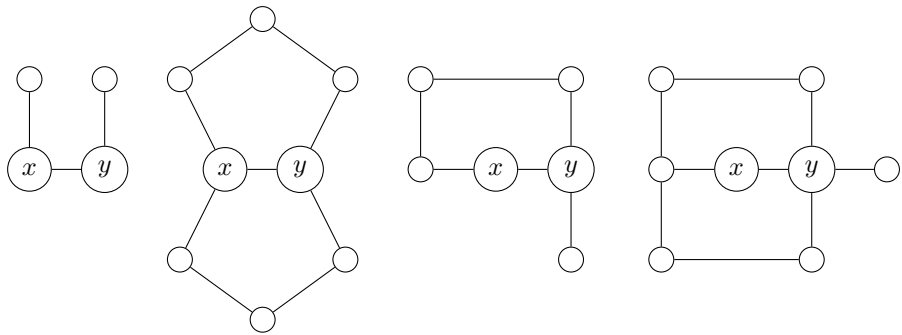
\begin{figure}[h]
\begin{center}
\begin{tikzpicture}[scale=0.2]

\node (vb) at (0,15) [circle,draw] {$y$};
\node (va) at (-5,15) [circle,draw] {$x$};
\node (vc) at (-5,21) [circle,draw] {};
\node (vd) at (0,21) [circle,draw] {};

\draw (va) -- (vb);
\draw (va) -- (vc);
\draw (vb) -- (vd);

\node (v0) at (16,21) [circle,draw] {};
\node (v1) at (13,15) [circle,draw] {$y$};
\node (v2) at (8,15) [circle,draw] {$x$};
\node (v3) at (5,21) [circle,draw] {};
\node (v4) at (10.5,25) [circle,draw] {};
\node (v5) at (5,9) [circle,draw] {};
\node (v6) at (10.5,5) [circle,draw] {};
\node (v7) at (16,9) [circle,draw] {};

\node (v8) at (21,15) [circle,draw] {};
\node (v9) at (26,15) [circle,draw] {$x$};
\node (v10) at (31,15) [circle,draw] {$y$};
\node (v11) at (21,21) [circle,draw] {};
\node (v12) at (31,21) [circle,draw] {};
\node (v13) at (31,9) [circle,draw] {};

\node (v14) at (45-8,15) [circle,draw] {};
\node (v15) at (50-8,15) [circle,draw] {$x$};
\node (v16) at (55-8,15) [circle,draw] {$y$};
\node (v17) at (60-8,15) [circle,draw] {};
\node (v18) at (45-8,9) [circle,draw] {};
\node (v19) at (55-8,9) [circle,draw] {};
\node (v20) at (45-8,21) [circle,draw] {};
\node (v21) at (55-8,21) [circle,draw] {};

\draw (v1) -- (v0);
\draw (v2) -- (v1);
\draw (v3) -- (v2);
\draw (v4) -- (v0);
\draw (v4) -- (v3);
\draw (v5) -- (v2);
\draw (v6) -- (v5);
\draw (v7) -- (v1);
\draw (v7) -- (v6);

\draw (v9) -- (v8);
\draw (v10) -- (v9);
\draw (v11) -- (v8);
\draw (v12) -- (v10);
\draw (v12) -- (v11);
\draw (v13) -- (v10);
\draw (v15) -- (v14);
\draw (v16) -- (v15);
\draw (v17) -- (v16);
\draw (v18) -- (v14);
\draw (v19) -- (v16);
\draw (v19) -- (v18);
\draw (v20) -- (v14);
\draw (v21) -- (v16);
\draw (v21) -- (v20);

\end{tikzpicture}
\end{center}
\caption{Local Structures for girth at least five} \label{fig_local5}
\end{figure}

The above lemma lends us important ideas in analyzing the the local structure of Ricci-flat graphs of girth $4$. Actually, we will obtain the following very useful lemmas.

\begin{lem}\label{lem_bound4}
Let $xy$ be an edge of a graph $G$, and $xy$ is in exactly one $4$-cycle but is not in any $3$-cycle. Then $\kappa(x,y) \leqslant \frac{2}{d_x}+\frac{2}{d_y}-1$.
\end{lem}

\begin{pf} Since the edge $xy$ is in a unique $4$-cycle, let $z$ be the other neighbor of $x$ in this cycle.

Let
\begin{equation*}
f(u) = \begin{cases}
0 &\text{if } u\in N[x] \setminus \{y,z\},\\
2 &\text{if } u \in N(y) \setminus \{x\},\\
1 & \text{otherwise}.
\end{cases}
\end{equation*}
Obviously, $f$ is a $1$-Lipschitz function over graph $G$. By lemma \ref{lem_dual},

\begin{eqnarray*}
W(m_x^\alpha, m_y^\alpha)  & \geqslant & \sum_{u\in V} f(u)[m_y^\alpha(u)-m_x^\alpha(u)] \\
 & =  & (\alpha - \frac{1-\alpha}{d_x})+(0-\frac{1-\alpha}{d_x})
 +2(d_y-1)(\frac{1-\alpha}{d_y}-0) \\
  & = & (2 - \alpha)-(1-\alpha)(\frac{2}{d_x}+\frac{2}{d_y}).
\end{eqnarray*}

So

\begin{equation*}
\kappa(x,y)=\lim_{\alpha \rightarrow 1}\frac{1-W(m_x^\alpha, m_y^\alpha) }{1-\alpha} \leqslant \frac{2}{d_x}+\frac{2}{d_y}-1.
\end{equation*}

\hfill$\square$
\end{pf}

The next lemma characterizes the local structures for edges in a $4$-cycle of a Ricci-flat graph $G$ with girth $4$ and disjoint $4$-cycles.

\begin{lem}\label{lem_local4}
Suppose that $G$ is a graph with girth $4$, and the $4$-cycles of $G$ are mutually vertex-disjoint. Let $xy$ be an edge of $G$ in a $4$-cycle with Ricci curvature $\kappa(x,y)=0$. Without loss of generality, we assume $d_x \leqslant d_y$, then one of the following statements holds.
\begin{enumerate}
\item $d_x=2, d_y=4$, and $xy$ is not in any $5$-cycle.
\item $d_x=d_y=3$, and $xy$ is not in any $5$-cycle.
\item $d_x=3, d_y=4$. Let $x_1$ and $x_2$ be the two neighbors of $x$ besides $y$ with $x_1$ in the $4$-cycle, and let $y_1$ and $y_2$ be the two neighbors of $y$ not in the $4$-cycle, then either $d(x_1,y_1)=d(x_2,y_2)=2$ (Type A), or $d(x_2,y_1)=d(x_2,y_2)=2$ (Type B).
\item $d_x=d_y=4$. Let $x_1$ and $x_2$ be the two neighbors of $x$ not in the $4$-cycle, and let $y_1$ and $y_2$ be the two neighbors of $y$ not in the $4$-cycle, then $d(x_1,y_1)=d(x_2,y_2)=2$.
\end{enumerate}
\end{lem}

\begin{figure}[h]
\begin{center}
\begin{tikzpicture}[scale=0.2]

\node (v0) at (12,81) [circle,draw] {$y$};
\node (v1) at (12,74) [circle,draw] {};
\node (v2) at (5,74) [circle,draw] {};
\node (v3) at (5,81) [circle,draw] {$x$};
\node (v4) at (19,85) [circle,draw] {};
\node (v5) at (19,81) [circle,draw] {};

\node (v6) at (30,81) [circle,draw] {$x$};
\node (v7) at (30,74) [circle,draw] {};
\node (v8) at (37,74) [circle,draw] {};
\node (v9) at (37,81) [circle,draw] {$y$};
\node (v10) at (24,85) [circle,draw] {};
\node (v11) at (43,85) [circle,draw] {};

\node (v12) at (50,81) [circle,draw] {$x$};
\node (v13) at (57,81) [circle,draw] {$y$};
\node (v14) at (57,74) [circle,draw] {};
\node (v15) at (50,74) [circle,draw] {};
\node (v16) at (57,70) [circle,draw] {};
\node (v17) at (65,78) [circle,draw] {};
\node (v18) at (65,84) [circle,draw] {};
\node (v19) at (57,87) [circle,draw] {};
\node (v20) at (50,87) [circle,draw] {};

\node at (64,72) {Type A};

\node (v21) at (7,60) [circle,draw] {$x$};
\node (v22) at (14,60) [circle,draw] {$y$};
\node (v23) at (7,53) [circle,draw] {};
\node (v24) at (14,53) [circle,draw] {};
\node (v25) at (7,66) [circle,draw] {};
\node (v26) at (14,70) [circle,draw] {};
\node (v27) at (22,64) [circle,draw] {};
\node (v28) at (22,60) [circle,draw] {};
\node (v29) at (14,66) [circle,draw] {};

\node at (20,55) {Type B};

\node (v30) at (38,60) [circle,draw] {$x$};
\node (v31) at (45,60) [circle,draw] {$y$};
\node (v32) at (45,53) [circle,draw] {};
\node (v33) at (38,53) [circle,draw] {};
\node (v34) at (30,64) [circle,draw] {};
\node (v35) at (30,60) [circle,draw] {};
\node (v36) at (53,64) [circle,draw] {};
\node (v37) at (53,60) [circle,draw] {};
\node (v38) at (41.5,66) [circle,draw] {};
\node (v39) at (41.5,70) [circle,draw] {};

\draw (v1) -- (v0);
\draw (v2) -- (v1);
\draw (v3) -- (v0);
\draw (v3) -- (v2);
\draw (v4) -- (v0);
\draw (v5) -- (v0);
\draw (v7) -- (v6);
\draw (v8) -- (v7);
\draw (v9) -- (v6);
\draw (v9) -- (v8);
\draw (v10) -- (v6);
\draw (v11) -- (v9);
\draw (v13) -- (v12);
\draw (v14) -- (v13);
\draw (v15) -- (v12);
\draw (v15) -- (v14);
\draw (v16) -- (v15);
\draw (v17) -- (v13);
\draw (v17) -- (v16);
\draw (v18) -- (v13);
\draw (v19) -- (v18);
\draw (v20) -- (v12);
\draw (v20) -- (v19);
\draw (v22) -- (v21);
\draw (v23) -- (v21);
\draw (v24) -- (v22);
\draw (v24) -- (v23);
\draw (v25) -- (v21);
\draw (v26) -- (v25);
\draw (v27) -- (v22);
\draw (v27) -- (v26);
\draw (v28) -- (v22);
\draw (v29) -- (v25);
\draw (v29) -- (v28);
\draw (v31) -- (v30);
\draw (v32) -- (v31);
\draw (v33) -- (v30);
\draw (v33) -- (v32);
\draw (v34) -- (v30);
\draw (v35) -- (v30);
\draw (v36) -- (v31);
\draw (v37) -- (v31);
\draw (v38) -- (v35);
\draw (v38) -- (v37);
\draw (v39) -- (v34);
\draw (v39) -- (v36);

\end{tikzpicture}
\end{center}
\caption{Local Structures in $4$-cycle} \label{fig_local4}
\end{figure}
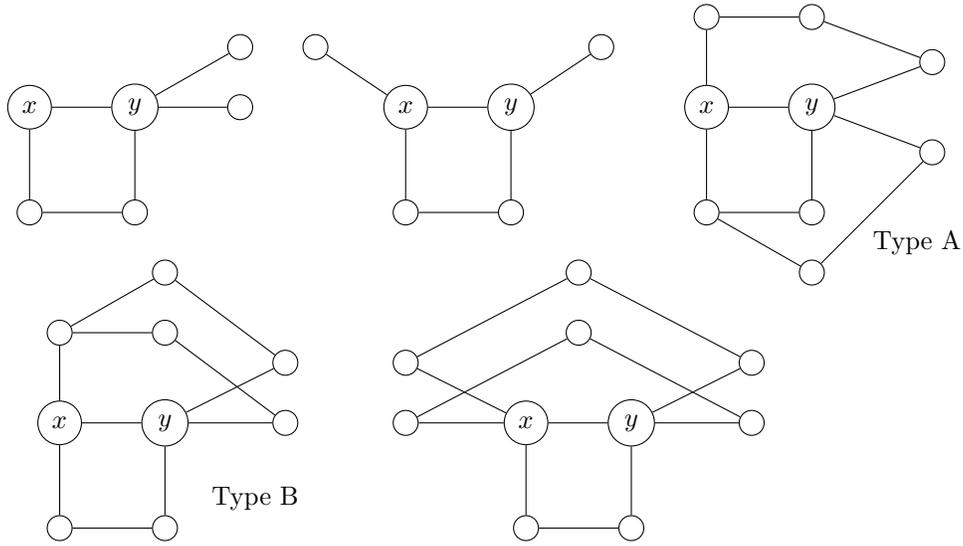

\begin{rmk}
Lemma \ref{lem_local4} claims that under certain conditions, for each edge $xy$ with Ricci curvature $0$, there are four possible degree combinations $\{d_x,d_y\}$. But there are two possible local structures in the $\{3,4\}$ combination, which are denoted by type A and type B, resulting in five local structures in all.
\end{rmk}

\begin{pf}
By Lemma \ref{lem_bound4} and the hypothesis that $\kappa(x,y)=0$, we have $0\leqslant \frac{2}{d_x}+\frac{2}{d_y}-1$. Solving this inequality, we have the following solutions.
\begin{enumerate}[label=(\Alph*)]
\item $d_x=2$, $d_y\geqslant 2$.
\item $d_x=3$, $d_y=3,4,5,6$.
\item $d_x=d_y=4$.
\end{enumerate}

Because the $4$-cycles of $G$ are disjoint, the vertex $y$ must be incident to another edge which is not in any $3$-cycles or $4$-cycles, if $d_y\geqslant 3$. By Lemma \ref{lem_local5}, $d_y\leqslant 4$. So the possible values of $d_x$ and $d_y$ become the following.

\begin{enumerate}[label=(\Alph*)]
\item $d_x=2$, $d_y=2,3,4$.
\item $d_x=3$, $d_y=3,4$.
\item $d_x=d_y=4$.
\end{enumerate}

Simple calculation shows that if $d_x=d_y=2$, then $\kappa(x,y)=1$. So no local structure is possible for this degree combination.

If $d_x=2$ and $d_y=3$, let $x_1$ be the other neighbor of $x$ besides $y$, and let $y_1$ be the neighbor of $y$ that is not in the $4$-cycle. Because the $4$-cycles of $G$ are disjiont, $d(x_1,y_1)\geqslant 2$. If $d(x_1,y_1)=2$, we have $\kappa(x,y)=\frac{1}{2}$. If $d(x_1,y_1)\geqslant 3$, we have $\kappa(x,y)=\frac{1}{3}$.

If $d_x=2$ and $d_y=4$, let $x_1$ be the other neighbor of $x$ besides $y$, and let $y_1$ and $y_2$ be the two neighbors of $y$ that is not in the $4$-cycle. If $d(x_1,y_1)=2$ or $d(x_1,y_2)=2$, then $\kappa(x,y)=\frac{1}{4}$. If $d(x_1,y_1)\geqslant 3$ and $d(x_1,y_2)\geqslant 3$, then $\kappa(x,y)=0$. Therefore, the edge $xy$ is not in any $5$-cycle.

The above calculations for the case (A) $d_x=2$ and $d_y=2,3,4$ can be summarized in Table \ref{tab_2x}.

\begin{table}
\begin{center}
  \begin{tabular}{ | c | c | c | c | c | }
    \hline
   $d_x$ & $d_y$ & $d(x_1,y_1)$ & $d(x_1,y_1),d(x_1,y_2)$ & $\kappa$ \\ \hline
   $2$ & $2$ & - & - & $1$ \\ \hline
   $2$ & $3$ & $2$ & - & $\frac{1}{2}$ \\ \hline
   $2$ & $3$ & $\geqslant 3$ & - & $\frac{1}{3}$ \\ \hline
   $2$ & $4$ & - & $2,\geqslant 2$ or $\geqslant 2, 2$ & $\frac{1}{4}$ \\ \hline
   $2$ & $4$ & - & $\geqslant 3,\geqslant 3$ & $0$ \\ \hline
  \end{tabular}
\end{center}
\caption{$d_x=2$}\label{tab_2x}
\end{table}

If $d_x=3$ and $d_y=3$, let $x_1$ and $x_2$ be the other neighbors of $x$ besides $y$, with $x_1$ in the $4$-cycle.
And let $y_1$ and $y_2$ be the two neighbors of $y$ besides $x$, with $y_1$ in the $4$-cycle. If $d(x_2,y_2)=2$, then $\kappa(x,y)=\frac{1}{3}$. If $d(x_2,y_2)\geqslant 3$, then $\kappa(x,y)=0$.

If $d_x=3$ and $d_y=4$,  let $x_1$ and $x_2$ be the other neighbors of $x$ besides $y$, with $x_1$ in the $4$-cycle.
And let $y_1$ and $y_2$ be the two neighbors of $y$ not in the $4$-cycle. Note that either $d(x_i,y_j)=2$  or  $d(x_i,y_j)=3$ for all $i,j=1,2$, so the complete calculations are divided into $8$ subcases according to the distances between $x_1, x_2$ and $y_1, y_2$. The result are listed in Table \ref{tab_34}. In three of the subcases (Lines 4,6,8 of the table), the Ricci curvature of edge $xy$ vanishes. Line 4 is Type B. Since the vertices $y_1$ and $y_2$ are interchangeable, line 6 and line 8 of the table can be combined to obtain the local structure of Type A.

\begin{table}

\begin{center}
  \begin{tabular}{ | c | c | c | }
    \hline
    $d(x_1,y_1),d(x_1,y_2)$ & $d(x_2,y_1),d(x_2,y_2)$ & $\kappa$ \\ \hline
    $3,3$ & $3,3$ & $-\frac{1}{3}$ \\ \hline
    $3,3$ & $2,3$ & $-\frac{1}{12}$ \\ \hline
    $3,3$ & $2,2$ & $0$ \\ \hline
    $2,\geqslant 2$ & $3,3$ & $-\frac{1}{4}$ \\ \hline
    $2,2$ & $2,3$ & $0$ \\ \hline
    $2,3$ & $2,3$ & $-\frac{1}{12}$ \\ \hline
    $2,\geqslant 2$ & $3,2$ & $0$ \\ \hline
    $2,\geqslant 2$ & $2,2$ & $\frac{1}{12}$ \\
    \hline
  \end{tabular}
\end{center}

\caption{$d_x=3, d_y=4$}\label{tab_34}
\end{table}

If $d_x=4$ and $d_y=4$,  let $x_1$ and $x_2$ be the two neighbors of $x$ not in the $4$-cycle.
And let $y_1$ and $y_2$ be the two neighbors of $y$ not in the $4$-cycle. Table \ref{tab_44} shows the Ricci curvature of edge $xy$ for different subcases. The unique subcase that the Ricci curvature vanishes is illustrated in bottom right of Figure \ref{fig_local4}.

\begin{table}
\begin{center}
  \begin{tabular}{ | c | c | c | }
    \hline
    $d(x_1,y_1),d(x_1,y_2)$ & $d(x_2,y_1),d(x_2,y_2)$ & $\kappa$ \\ \hline
    $3,3$ & $3,3$ & $-\frac{1}{2}$ \\ \hline
    $2,3$ & $\geqslant 2,3$ & $-\frac{1}{4}$ \\ \hline
    $2,\geqslant 2$ & $\geqslant 2,2$ & $0$ \\
    \hline
  \end{tabular}
\end{center}
\caption{$d_x=4, d_y=4$}\label{tab_44}
\end{table}

\hfill$\square$

\end{pf}

Lemma \ref{lem_local5} and Lemma \ref{lem_local4} will be applied repeatedly in proving the main result in the next section.


\section{The Main Result}\label{sec_main}

This section proves Theorem \ref{major} by exhausting all possible cases.
\begin{pf2}We start by investigating a $4$-cycle of $G$. By Lemma \ref{lem_local4}, the degree sequence of a $4$-cycle of $G$ in cyclic order can be only one of the following cases.

\begin{enumerate}
\item (2,4,2,4)
\item (2,4,4,4)
\item (3,3,3,3)
\item (3,3,3,4)
\item (3,3,4,4)
\item (3,4,4,4)
\item (3,4,3,4)
\item (4,4,4,4)
\end{enumerate}

We will show that in the first six cases, the graph $G$ could not exist. And in the last two cases, exactly one graph is possible for each case.

\textbf{Case 1. (2,4,2,4)}. Let $a,b,c,d$ be the four vertices of the $4$-cycle, in the order of the degree sequence. That is $d(a)=d(c)=2$ and $d(b)=d(d)=4$. Let $b_1$ and $b_2$ the other two neighbors of $b$, and let $d_1$ and $d_2$ be the other two neighbors of $d$. Obviously, $b_i$ and $d_j$ ($1\leqslant i,j \leqslant 2$) are distinct vertices, otherwise there would be $4$-cycles with common edges. In the remaining cases, we will denote and refer to the vertices in the $4$-cycle and their neighbors in a similar manner.

Because the edge $bb_1$ does not lie in any $4$-cycle by the hypothesis of the theorem, so it must satisfy the local structure of Lemma \ref{lem_local5}. Since $d(b)=4$, so $d(b_1)=2$. By the same reason $d(b_2)=d(d_1)=d(d_2)=2$. Let $z$ be the other neighbor of $b_1$ besides $b$. See Figure \ref{fig_case1}. Note that $z$ must be distinct from $d_1$ or $d_2$. Suppose to the contrast that the other neighbor of $b_1$ is $d_1$, then the edge $b_1d_1$ does not satisfy Lemma \ref{lem_local5}.

Now we apply Lemma \ref{lem_local5} to edge $bb_1$, at least two vertices of $a,c,b_2$ have distance $2$ from $z$. But this is impossible (because there is no way to form a $2$-path from $z$ to either $a$ or $c$), so no graph exists for this case.

\begin{figure}[h]
\begin{center}
\begin{tikzpicture}[scale=0.15]

\node (v0) at (30,30) [circle,draw] {$b$};
\node (v1) at (0,30) [circle,draw] {$a$};
\node (v2) at (0,0) [circle,draw] {$d$};
\node (v3) at (30,0) [circle,draw] {$c$};
\node (v4) at (8,5) [circle,draw] {$d_1$};

\node (v6) at (8,12) [circle,draw] {$d_2$};

\node (v8) at (15,15) [circle,draw] {$z$};

\node (v10) at (22,25) [circle,draw] {$b_2$};
\node (v11) at (22,18) [circle,draw] {$b_1$};

\draw (v1) -- (v0);
\draw (v2) -- (v1);
\draw (v3) -- (v0);
\draw (v3) -- (v2);

\draw (v4) -- (v2);
\draw (v6) -- (v2);
\draw (v10) -- (v0);
\draw (v11) -- (v0);
\draw (v11) -- (v8);

\end{tikzpicture}
\end{center}
\caption{Case 1. (2,4,2,4)} \label{fig_case1}
\end{figure}
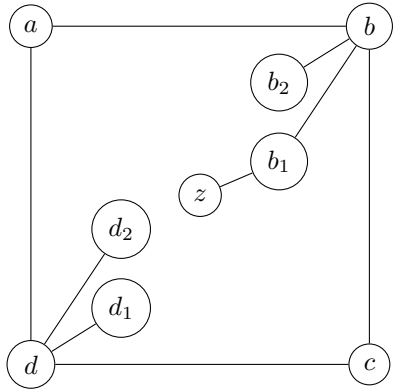

\textbf{Case 2. (2,4,4,4)}. See Figure \ref{fig_case2}. The same as Case 1, because the degree of vertices $b,c,d$ are $4$, the degree of vertices $b_i,c_i,d_i$ ($i=1,2$) are all $2$. By applying Lemma \ref{lem_local4} to edge $bc$, without loss of generality, let $z_i$ be the common neighbor of $b_i$ and $c_i$ ($i=1,2$). By applying Lemma \ref{lem_local4} again to edge $cd$, we know that the vertices $c_i$ and $d_i$ have a common neighbor, for $i=1,2$. But since all the vertices $c_1,c_2,d_1,d_2$ have degree $2$, the common neighbor of $c_i$ and $d_i$ has to be $z_i$, for $i=1,2$. Now, all vertices in Figure \ref{fig_case2} cannot be extended except $z_1$ and $z_2$. Therefore, the edge $bb_1$ does not satisfy Lemma \ref{lem_local5} (the edge $bb_1$ does not lie in two $5$-cycles), no graph exists for this case, either.

\begin{figure}[h]
\begin{center}
\begin{tikzpicture}[scale=0.15]

\node (v0) at (30,30) [circle,draw] {$b$};
\node (v1) at (0,30) [circle,draw] {$a$};
\node (v2) at (0,0) [circle,draw] {$d$};
\node (v3) at (30,0) [circle,draw] {$c$};
\node (v4) at (8,5) [circle,draw] {$d_1$};
\node (v5) at (22,5) [circle,draw] {$c_1$};
\node (v6) at (8,12) [circle,draw] {$d_2$};
\node (v7) at (22,12) [circle,draw] {$c_2$};
\node (v8) at (15,10) [circle,draw] {$z_1$};
\node (v9) at (15,20) [circle,draw] {$z_2$};
\node (v10) at (22,25) [circle,draw] {$b_2$};
\node (v11) at (22,18) [circle,draw] {$b_1$};

\draw (v1) -- (v0);
\draw (v2) -- (v1);
\draw (v3) -- (v0);
\draw (v3) -- (v2);
\draw (v4) -- (v2);
\draw (v5) -- (v3);
\draw (v6) -- (v2);
\draw (v7) -- (v3);
\draw (v8) -- (v4);
\draw (v8) -- (v5);
\draw (v9) -- (v6);
\draw (v9) -- (v7);
\draw (v10) -- (v0);
\draw (v10) -- (v9);
\draw (v11) -- (v0);
\draw (v11) -- (v8);

\end{tikzpicture}
\end{center}
\caption{Case 2. (2,4,4,4)} \label{fig_case2}
\end{figure}
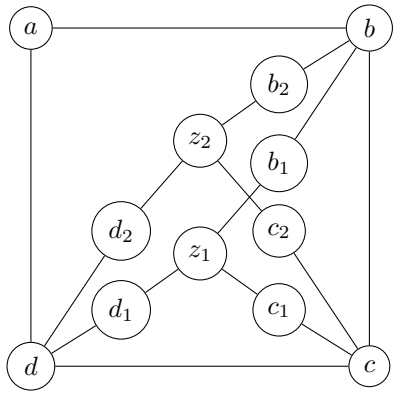

\textbf{Case 3. (3,3,3,3)}. In this case, each vertex of the $4$-cycle, $a, b, c$ and $d$, has exactly one neighbor outside the $4$-cycle, denoted by $a_1, b_1, c_1$ and $d_1$, respectively. By Lemma \ref{lem_local5}, the degree of $a_1, b_1, c_1$ and $d_1$ can be either $2$ or $3$.

If $d(b_1)=3$, by Lemma \ref{lem_local5}, the edge $bb_1$ is shared by two $5$-cycles, this contradicts with the fact the the edge $ab$ cannot lie in any $5$-cycles. If $d(b_1)=2$, by Lemma \ref{lem_local5}, the edge $bb_1$ need to form a $5$-cycle with either $ab$ or $bc$, which is also a contradiction. So no graph exists for this case.

\textbf{Case 4. (3,3,3,4)}. The same as Case 3, by applying Lemma \ref{lem_local5} to edge $bb_1$, there will be a contradiction. So no graphs exists for this case.

\textbf{Case 5. (3,3,4,4)}. See Figure \ref{fig_case5}. Both $a_1$ and $b_1$ must have degree $2$, otherwise by the same argument in Case 3 the edge $ab$ lies in a $5$-cycle, a contradiction. Also, the degree of $c_1, c_2, d_1$ and $d_2$ are all $2$. By applying Lemma \ref{lem_local4} to edge $cd$, let $z_i$ be the common neighbor of $c_i$ and $d_i$, for $i=1,2$. Thus the edges $bc$ and $da$ have no way to satisfy the local condition. Again, no graph exists for this case.

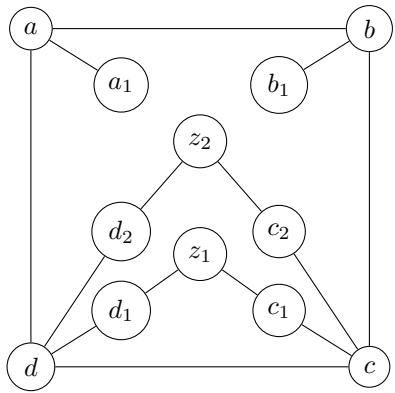
\begin{figure}[h]
\begin{center}
\begin{tikzpicture}[scale=0.15]

\node (v0) at (30,30) [circle,draw] {$b$};
\node (v1) at (0,30) [circle,draw] {$a$};
\node (v2) at (0,0) [circle,draw] {$d$};
\node (v3) at (30,0) [circle,draw] {$c$};
\node (v4) at (8,5) [circle,draw] {$d_1$};
\node (v5) at (22,5) [circle,draw] {$c_1$};
\node (v6) at (8,12) [circle,draw] {$d_2$};
\node (v7) at (22,12) [circle,draw] {$c_2$};
\node (v8) at (15,10) [circle,draw] {$z_1$};
\node (v9) at (15,20) [circle,draw] {$z_2$};
\node (v10) at (22,25) [circle,draw] {$b_1$};
\node (v11) at (8,25) [circle,draw] {$a_1$};

\draw (v1) -- (v0);
\draw (v2) -- (v1);
\draw (v3) -- (v0);
\draw (v3) -- (v2);
\draw (v4) -- (v2);
\draw (v5) -- (v3);
\draw (v6) -- (v2);
\draw (v7) -- (v3);
\draw (v8) -- (v4);
\draw (v8) -- (v5);
\draw (v9) -- (v6);
\draw (v9) -- (v7);
\draw (v10) -- (v0);
\draw (v11) -- (v1);

\end{tikzpicture}
\end{center}
\caption{Case 5. (3,3,4,4)} \label{fig_case5}
\end{figure}

\textbf{Case 6. (3,4,4,4)}. The structure of the graph is similar to that of Case 2, except that $a$ will have a neighbor $a_1$. To satisfy the local condition for edges $ab$ and $da$, the vertex $a_1$ must be adjacent to both $z_1$ and $z_2$, see Figure \ref{fig_case6}. But then the edge $a_1z_1$ does not satisfy the local condition. So no graph exists for this case.

\begin{figure}[h]
\begin{center}
\begin{tikzpicture}[scale=0.15]

\node (v0) at (30,30) [circle,draw] {$b$};
\node (v1) at (0,30) [circle,draw] {$a$};
\node (v2) at (0,0) [circle,draw] {$d$};
\node (v3) at (30,0) [circle,draw] {$c$};
\node (v4) at (8,5) [circle,draw] {$d_1$};
\node (v5) at (22,5) [circle,draw] {$c_1$};
\node (v6) at (8,12) [circle,draw] {$d_2$};
\node (v7) at (22,12) [circle,draw] {$c_2$};
\node (v8) at (15,10) [circle,draw] {$z_1$};
\node (v9) at (15,20) [circle,draw] {$z_2$};
\node (v10) at (22,25) [circle,draw] {$b_2$};
\node (v11) at (22,18) [circle,draw] {$b_1$};
\node (v12) at (8,25) [circle,draw] {$a_1$};

\draw (v1) -- (v12);
\draw (v8) -- (v12);
\draw (v9) -- (v12);
\draw (v1) -- (v0);
\draw (v2) -- (v1);
\draw (v3) -- (v0);
\draw (v3) -- (v2);
\draw (v4) -- (v2);
\draw (v5) -- (v3);
\draw (v6) -- (v2);
\draw (v7) -- (v3);
\draw (v8) -- (v4);
\draw (v8) -- (v5);
\draw (v9) -- (v6);
\draw (v9) -- (v7);
\draw (v10) -- (v0);
\draw (v10) -- (v9);
\draw (v11) -- (v0);
\draw (v11) -- (v8);

\end{tikzpicture}
\end{center}
\caption{Case 6. (3,4,4,4)} \label{fig_case6}
\end{figure}
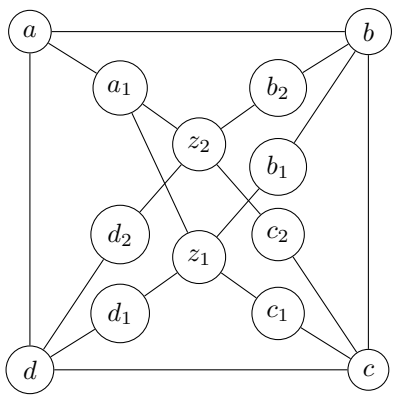

\textbf{Case 7. (3,4,3,4)}. The degree combination for the two vertices of each edge in the $4$-cycle is $\{3,4\}$. There are two types of local structures for the $\{3,4\}$ combination, namely type A and type B. It is easy to show that the four edges in the $4$-cycle must satisfy the same type of local condition. If all of them are type A, no graph is possible. If all of them are type B, we obtain the graph $R_2$, see Figure \ref{fig_main1b}.

\textbf{Case 8. (4,4,4,4)}. By applying Lemma \ref{lem_local4} to edge $ab$, let $z_i$ be the common vertex of $a_i$ and $b_i$ ($i=1,2$), respectively (See Figure \ref{fig_case8}). Then by applying Lemma \ref{lem_local5} to edge $bb_1$, the degree of $b_1$ must be two. In other words, $b_1$ has no other neighbors besides $b$ and $z_1$. Now by applying Lemma \ref{lem_local4} to edges $bc$ and $cd$, we obtain the graph $R_1$, see Figure \ref{fig_main1b}.

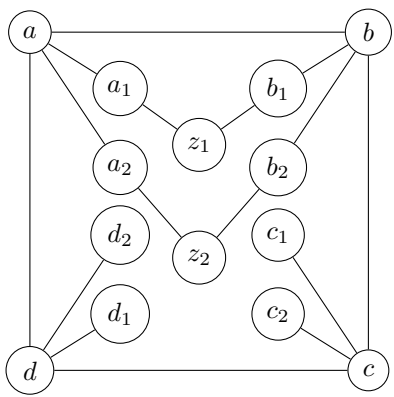
\begin{figure}[h]
\begin{center}
\begin{tikzpicture}[scale=0.15]

\node (v0) at (30,30) [circle,draw] {$b$};
\node (v1) at (0,30) [circle,draw] {$a$};
\node (v2) at (0,0) [circle,draw] {$d$};
\node (v3) at (30,0) [circle,draw] {$c$};
\node (v5) at (22,5) [circle,draw] {$c_2$};
\node (v6) at (8,18) [circle,draw] {$a_2$};
\node (v7) at (22,12) [circle,draw] {$c_1$};
\node (v8) at (15,10) [circle,draw] {$z_2$};
\node (v9) at (15,20) [circle,draw] {$z_1$};
\node (v10) at (22,25) [circle,draw] {$b_1$};
\node (v11) at (22,18) [circle,draw] {$b_2$};
\node (v12) at (8,25) [circle,draw] {$a_1$};

\node (v4) at (8,5) [circle,draw] {$d_1$};
\node (v13) at (8,12) [circle,draw] {$d_2$};

\draw (v2) -- (v4);
\draw (v2) -- (v13);

\draw (v1) -- (v12);

\draw (v9) -- (v12);
\draw (v1) -- (v0);
\draw (v2) -- (v1);
\draw (v3) -- (v0);
\draw (v3) -- (v2);

\draw (v1) -- (v6);
\draw (v8) -- (v6);
\draw (v5) -- (v3);
\draw (v7) -- (v3);

\draw (v10) -- (v0);
\draw (v10) -- (v9);
\draw (v11) -- (v0);
\draw (v11) -- (v8);

\end{tikzpicture}
\end{center}
\caption{Case 8. (8,4,4,4)} \label{fig_case8}
\end{figure}

Finally, it is easy to check that the graphs $R_1$ and $R_2$ are indeed Ricci-flat.
\hfill$\square$
\end{pf2}

\begin{rmk}[\bf Further Study]
Our method in proving Theorem \ref{major} might be extended further to study the Ricci-flat graphs of girth $4$ and with edge-disjoint $4$-cycles. In such an extension, more involved discussions are expected, especially when one wants to generalize the result of Lemma \ref{lem_local4} so that more local degree combinations for an edge $xy$ in a $4$-cycle are included. By the examples given in Figure \ref{fig_edge}, we see that any characterization of Ricci-flat graphs of girth $4$ and with edge-disjoint $4$-cycles must contain infinitely many types.
\end{rmk}

\section*{References}

\bibliography{1}

\end{document}